\documentclass{amsart}
\usepackage{amsmath,amsthm,amssymb,IMjournal}
\usepackage{graphicx,subfig}
\usepackage{color}
\usepackage{float}
\usepackage{accents}
\usepackage{tikz}
\usetikzlibrary{intersections,calc}
\begin{document}
\newcommand{\R}{\mathbb R}
\newcommand{\N}{\mathbb N}
\newcommand{\Z}{\mathbb Z}
\newcommand{\dd}{\mathrm{d}}
\newcommand{\lam}{\lambda}
\newcommand{\bfx}{\mathbf{x}}
\newcommand{\bfy}{\mathbf{y}}
\newcommand{\bfs}{\mathbf{s}}
\newcommand{\bfv}{\mathbf{v}}
\newcommand{\bfz}{\mathbf{z}}
\newcommand{\bfa}{\mathbf{a}}
\newcommand{\bfb}{\mathbf{b}}
\newcommand{\bfq}{\mathbf{q}}
\newcommand{\hK}{\widehat{K}}
\newcommand{\tK}{\widetilde{K}}
\newcommand{\cK}{\Check{K}}
\newcommand{\wA}{\widetilde{A}}
\newcommand{\hv}{\hat{v}}
\newcommand{\T}{\mathcal{T}}
\newcommand{\I}{\mathcal{I}}
\newcommand{\OO}{\mathcal{O}}
\newcommand{\PP}{\mathcal{P}}
\newcommand{\NN}{\mathcal{N}}
\newcommand{\Kab}{K_{ab}}
\newcommand{\bK}{\mathbf{K}}
\newcommand{\ux}{\underbar{x}}
\newcommand{\ox}{\overline{\mathrm{x}}}
\newcommand{\vvskip}{\vspace{5pt}}
\newcommand{\gmax}{\gamma_{\max}}

\renewcommand{\figurename}{Figure~}

\newtheorem{theorem}{Theorem}
\newtheorem{lemma}[theorem]{Lemma}
\newtheorem{definition}[theorem]{Definition}
\newtheorem{corollary}[theorem]{Corollary}
\newtheorem{assumption}[theorem]{Assumption}
\newtheorem{condition}[theorem]{Condition}

\newtheorem{The}{Theorem}[section]

\numberwithin{equation}{section}

\title{Error estimation for finite element solutions on meshes 
that contain thin elements}

\author{\|Kenta |Kobayashi|, Kunitachi,
        \|Takuya |Tsuchiya|, Toyonaka}

\rec{Feburary 29, 2024}

\dedicatory{Cordially dedicated to the memory of I.$\,$Babu\v{s}ka}\enddedicatory

\abstract 
%%%%%%
%%%7 to 12 lines optimum; references in the abstract should be given in full form, e.g. V. Nov\'ak, M. Novotn\'y: Linear extensions of orderings. Czech. Math. J. 50 (2000), 853--864.
%%%%%
   In an error estimation of finite element solutions to the Poisson
   equation, we usually impose the shape regularity assumption on
   the meshes to be used.
   In this paper, we show that even if the shape regularity condition
   is violated, the standard error estimation can be obtained if ``bad''
   elements elements that violate the shape regularity or maximum angle
   condition are covered virtually by simplices {that satisfy the
   minimum angle condition}.
   A numerical experiment illustrates the theoretical result.
%   A numerical experiment confirms the theoretical result.
\endabstract



\thanks
   {}The authors were supported by JSPS KAKENHI Grant Numbers
   24K00538 and 21K03372. 
\endthanks

\section{Introduction}\label{intro}
The mathematical theory of finite element methods has been
well established, and many well written textbooks on the subject are
available \cite{BreSco, Cia78, EG}.
In these textbooks, the central assumption of the theory is the 
\textit{shape regularity condition} on the meshes that are used to
define finite element solutions.  The shape regularity condition
requires that each element is not very degenerate.  For example, if a
mesh consists of triangles, each triangle should not be very ``flat''.
It is known that the shape regularity condition is equivalent to the
\textit{minimum angle condition}.

However, suppose the problem domain has a highly complex shape, and/or
the problem has an anisotropic nature,  and we introduce an adaptive
mesh refinement procedure. In that case, the meshes might contain some
very thin elements and the error estimation of the finite element method
may thus lose its theoretical validity.

There are two ways to overcome this issue.  One way is to introduce
a new geometric condition to meshes, where the condition ensures a
standard convergence rate of finite element solutions.  It is well known
that if the maximum angles of triangles in a mesh are bounded by a fixed
constant $\theta < \pi$, then the standard error estimation of finite
element solutions holds \cite{BabAzi76}.  This condition is known as the
\textit{maximum angle condition}. The approach has been extended in many
ways.  See, for example,  \cite{Apel,Jamet, KobTsu16, Kri91}.

The other way was initiated by {Ku\v{c}era \cite{Kucera}}
and Duprez--Lleras--Lozinski \cite{DLL}.
Suppose that a mesh contains some ``bad'' elements. {In 
\cite{Kucera, DLL}, if bad elements satisfy some conditions, the usual
Lagrange interpolation is modified in certain ways on bad elements so
that the standard error estimation holds.} 

This paper extends their approach.  We show that if all
bad elements are covered virtually by simplices 
{that satisfy the minimum angle condition}, the
standard error estimation holds.  Note that the simplices are
virtual and do not need to be elements of the mesh considered.

In Section~2, we present the model problem, the precise definitions of
conformal meshes, and the above-mentioned conditions.  In addition, we
give an example of bad meshes for which the standard error estimation
holds.  Section~3 presents detailed assumptions of meshes and
the main theorem.   
Section~4 presents a simpler version of the assumption and the main
theorem for the case $n = 2$.
{In Section~5, we compare our approach with approaches adopted in
 prior works \cite{Kucera,DLL}.}
In Section~6, we present a numerical
experiment that confirms the theoretical result. 

\section{Preliminary} \label{preli}
\subsection{Model problem}
Let $\Omega \subset \R^n$ $(n=2,3)$ be a {bounded}
polygonal or polyhedral domain. For simplicity, we consider
the Poisson equation 
\begin{align*}
  - \Delta u = f \; \text{ in } \;\Omega, \qquad
     u = 0 \; \text{ on } \; \partial\Omega
\end{align*}
{and} its weak form
\begin{align}
   \int_\Omega \nabla u\cdot\nabla v \,\dd \bfx = 
     \int_\Omega fv \, \dd \bfx \quad
       \forall v \in H_0^1(\Omega)
     \label{exact-sol}
\end{align}
as the model problem, where $f \in L^2(\Omega)$ is {a} given function.
Here, $L^p(\Omega)$, $H^m(\Omega)$, $H_0^1(\Omega)$, and
$W^{m,p}(\Omega)$ ($m=1,2$, $1 \le p \le \infty$)  are
the standard Lebesgue and Sobolev spaces with norms and seminorms
$\|\cdot\|_{L^p(\Omega)}$, $\|\cdot\|_{H^m(\Omega)}$,
$|\cdot|_{H^m(\Omega)}$, $\|\cdot\|_{W^{m,p}(\Omega)}$, and 
$|\cdot|_{W^{m,p}(\Omega)}$.
In this paper, we assume that the exact solution $u$ of
\eqref{exact-sol} belongs to $H^2(\Omega)$.

\subsection{Finite element method}
We approximate the exact solution of \eqref{exact-sol} using the
standard piecewise linear finite element method.  To this end, we 
define the piecewise linear finite element spaces by
\begin{align*}
  S_h := \left\{v_h \in C^0(\overline{\Omega}) \bigm|
     v_h|_\tau \in \mathcal{P}_1, \forall \tau \in \T_h\right\},
 \quad
  S_{h0} := S_h \cap H_0^1(\Omega),
\end{align*}
where $\PP_1$ is the set of polynomials with degree at most $1$.
Here, $\T_h$ is a conformal simplicial mesh of $\Omega$.
That is, $\T_h$ is a finite 
set of $n$-simplices that satisfies the conditions
\begin{itemize}
\item $\displaystyle \bigcup_{\tau \in \T_h} \tau = \overline{\Omega}$.
\item If $\tau_1 \cap \tau_2 \neq \emptyset$ for
$\tau_1$, $\tau_2 \in \T_h$ with $\tau_1 \neq \tau_2$, then
    $\tau_1 \cap \tau_2$ is their common $r$-face ($r < n$)
    (that is, there are no hanging nodes in $\T_h$).
\end{itemize}
As usual, we define $h_\tau := \mathrm{diam}\,\tau$ and
$h := \max_{\tau \in \T_h} h_\tau$. 
The finite element solution $u_h \in S_{h0}$ is then defined as the
unique solution of the equation
\begin{align}
  \int_\Omega \nabla u_h\cdot \nabla v_h \, \dd \bfx =
  \int_\Omega f v_h\, \dd \bfx, \quad
    \forall v_h \in S_{h0}. 
   \label{FEsol}
\end{align}

\subsection{C\'ea's lemma and interpolations}
For the error of $u_h$, C\'ea's lemma claims that
{there exists a constant $C > 0$ such that}
\begin{align*}
  |u -  u_h|_{H^1(\Omega)} \le {C} \inf_{w_h \in S_{h0}}
     |u - w_ h|_{H^1(\Omega)}
     \le {C}\, |u - \Pi^*u|_{H^1(\Omega)},
\end{align*}
where $\Pi^*u \in S_h$ is \textit{any} interpolation of the
exact solution $u$.  Therefore, if we have
\begin{align}
  |u - \Pi^* u|_{H^1(\Omega)} \le C h |u|_{H^2(\Omega)},
  \label{inter}
\end{align}
we obtain the standard convergence rate
\begin{align*}
   |u - u_h|_{H^1(\Omega)} \le C h |u|_{H^2(\Omega)}.
\end{align*}

Usually, the Lagrange interpolation {operator} is used as $\Pi^*$.
The piecewise linear Lagrange interpolation $\Pi_\tau^1 u \in \PP_1$
is defined by
\begin{align*}
  (\Pi_\tau^1u)(\bfx_i) & = u(\bfx_i) \text{ on $n$-simplex } \tau
    \text{ with vertices } \bfx_i, i = 1,\dots,n+1, \\
  \Pi_h^1 u\bigm|_\tau & = \Pi_\tau^1 u, \qquad \forall \tau \in \T_h.
\end{align*}
In the next subsection, we mention conditions for which the
Lagrange interpolation {satisfies} the standard error estimation
\begin{align}
   |u - \Pi_\tau^1 u|_{W^{1,p}(\tau)} \le C_p h_\tau |u|_{W^{2,p}(\tau)},
   \quad \forall u \in W^{2,p}(\tau).
   \label{Lag-error}
\end{align}

\subsection{Conditions and constants related to meshes}

\begin{definition}
For an $n$-simplex $\tau$, 
%we define $\theta_\tau^{min}$ and
%$\theta_\tau^{max}$ as follows.
$\theta_\tau^{min}$ and $\theta_\tau^{max}$ denote the minimum
and maximum inner angles of $\tau$, respectively.
Here, the inner angles of a tetrahedron mean,
for $\theta_\tau^{min}$, the inner angles of each face and the solid
angles at vertices, and for $\theta_\tau^{max}$,
the inner angles of each face and the dihedral angles
\cite{BraKorKri08, Kri92}.
\end{definition}

The following conditions are well known for geometric conditions of
simplices.
\begin{definition}
$(1)$  If $\theta_\tau^{min} \ge \alpha > 0$,
$\tau$ is said to satisfy the \textbf{minimum angle condition} with
respect to $\alpha$.  \\
$(2)$ If $\theta_\tau^{max} \le \beta < \pi$,
$\tau$ is said to satisfy the \textbf{maximum angle condition} with
respect to $\beta$.
\end{definition}

If $\tau$ satisfies the minimum angle condition
with respect to $\alpha$,
the error estimation \eqref{Lag-error} holds for any $p$,
$1 \le p \le \infty$ with $C_p = C_{\alpha,p}^{min}$.
Moreover, there exists a constant $D_\alpha$ such that
\begin{align}
   \frac{h_\tau}{\rho_\tau} \le D_\alpha,
   \label{SR}
\end{align}
where $\rho_\tau$ is the diameter of the inscribed ball in $\tau$.
This condition is known as the \textbf{shape regularity condition} with
respect to $D_\alpha$.   It is well known that the shape regularity
condition is equivalent to the minimum angle condition.

Suppose that $\tau$ satisfies the maximum angle condition
with respect to $\beta$.
If $n=2$, \eqref{Lag-error} holds for any $p$,
$1 \le p \le \infty$ with $C_p = C_{\beta,p}^{max}$.
If $n=3$, \eqref{Lag-error} holds for $p$,
$2 < p \le \infty$ with $C_p = C_{\beta,p}^{max}$; however,
$C_{\beta,p}^{max} \to \infty$ as $p \searrow 2${, see}
\cite{Duran, KobTsu20, Kri92, She94}.

To specify the characteristics of meshes, we introduce constants
$\theta$, $\psi$, where $0 < \theta \le \pi/3$, $0 < \psi \le \pi/3$,
and fix them throughout this paper.  We define
\begin{align*}
  K_\theta &:= 
  \begin{cases} \;
   \{\tau \subset \R^2 \mid \tau \text{ is a triangle and }
   \theta_\tau^{max} \le \pi - {2}\theta\},
    & \text{for $n = 2$},  \\
   \; \{\tau \subset \R^3 \mid 
    \tau \text{ is a tetrahedron and }
   \theta_\tau^{min} \ge \theta\},
   &  \text{for $n = 3$}.
   \end{cases}
\end{align*}
{A simplex $\tau$ that does not belong to $K_\theta$ is
regarded as ``bad'', a simplex $\tau$ that satisfies the minimum
angle condition $\theta_\tau^{min} \ge \theta$ is regarded  as ``good'',
and other types of simplex are regarded as ``ordinary''. 
Note that for $n = 2$, $K_\theta$ consists of ``good'' and ``ordinary''
triangles, and for $n = 3$,
all tetrahedrons in $K_\theta$ are ``good''.  }

For a simplex $\tau$ that satisfies $\tau \in K_\theta$ or
$\theta_\tau^{min} \ge \psi$, 
we rewrite \eqref{Lag-error} as
\begin{align}
  |u - \Pi_\tau^1 u|_{H^1(\tau)}
  \le A_{\theta,\psi} h_\tau|u|_{H^2(\tau)},
  \qquad \forall u \in H^2(\tau),
    \label{eq2.6}
\end{align}
where $A_{\theta,\psi}$ is a constant that depends only on $\theta$
and $\psi$.  {The existence of $A_{\theta,\psi}$ is assured by
the minimum and maximum angle conditions.} 
{
(The parameter $\theta$ is used to categorize elements in mesh
whereas the parameter $\psi$ relates to virtual simplices
that appear later.)
}

\subsection{Sobolev embedding theorem}
Let $\hat{\tau}$ be a reference simplex with $h_{\hat{\tau}} = 1$.
Then, by the Sobolev embedding
theorem \cite{Brezis}, we have
\begin{align*}
    H^2(\hat{\tau}) \subset W^{1,6}(\hat{\tau})
    \subset L^\infty(\hat{\tau})
\end{align*}
and
\begin{align*}
  \|\hat{u} - \Pi_{\hat{\tau}}^1 \hat{u}\|_{L^\infty(\hat{\tau})} 
   \le \widehat{B} |\hat{u} - \Pi_{\hat{\tau}}^1
    \hat{u}|_{H^2(\hat{\tau})}
    = \widehat{B} |\hat{u}|_{H^2(\hat{\tau})}, \quad
   \hat{u} \in H^2(\hat{\tau}),
\end{align*}
where $\widehat{B}$ is a constant that depends only on $\hat{\tau}$.
Here, we used the fact that $(u - \Pi_{\hat{\tau}} u)(\bfx_i) = 0$,
$i = 1, \dots, n+1$, and 
$\|\hat{u} - \Pi_{\hat{\tau}}^1\hat{u}\|_{H^2(\hat{\tau})}$ may be
bounded by 
$|\hat{u} - \Pi_{\hat{\tau}}^1\hat{u}|_{H^2(\hat{\tau})}$
up to a constant.
For a general simplex $\tau$ with $\theta_\tau^{min} \ge \psi$,
we consider an affine map $\eta: \hat{\tau} \to \tau$, 
$\eta(\bfx) = A\bfx + \bfb$, where
$A$ is an $n\times n$ {nonsingular} matrix and $\bfb \in \R^n$.
Each $u \in H^2(\tau)$ is pulled back to $\hat{u}$ as
$\hat{u}(\bfx) := u(\eta(\bfx))$.  We then have
\begin{align*}
   \|u - \Pi_\tau^1 u\|_{L^\infty(\tau)}^2 & =
   \|\hat{u} - \Pi_{\hat{\tau}}^1 \hat{u}\|_{L^\infty(\hat{\tau})}^2 
   \le \widehat{B}^2 |\hat{u}|_{H^2(\hat{u})}^2  \\
 & \le \widehat{B}^2 |\det A|^{-1}\|A\|^4  |u|_{H^2(\tau)}^2,
\end{align*}
where $\|A\|$ is the matrix norm of $A$ with respect to the Euclidean
norm.   We denote the Lebesgue measure of $\tau$ by $|\tau|$. By
\cite{Cia78}, the coefficient of the last term is estimated as
\begin{align*}
   \det A = \frac{|\tau|}{|\hat{\tau}|} \ge C h^n, \qquad
   \|A\| \le \frac{h_\tau}{\rho_{\hat{\tau}}},
\end{align*}
and we finally obtain the embedding inequality
\begin{align}
  \|u - \Pi_\tau^1 u\|_{L^\infty(\tau)}
  \le B_\psi h_\tau^{2-n/2}  |u|_{H^2(\tau)}, \quad
   \forall u \in H^2(\tau),
  \label{eq2.8}
\end{align}
where $B_\psi$ is a constant that depends on $\psi$.

\subsection{Example of error estimation on a bad mesh.}
\begin{figure}[htbp]
\centering
\begin{tikzpicture}[line width = 0.8pt,scale=0.5]
   \coordinate (A1) at (0.0,0.0);
   \coordinate (A2) at (2.0,0.0);
   \coordinate (A3) at (4.0,0.0);
   \coordinate (A4) at (6.0,0.0);
   \coordinate (A5) at (8.0,0.0);
   \coordinate (B1) at (0.0,2.0);
   \coordinate (B5) at (8.0,2.0);
   \coordinate (C1) at (0.0,4.0);
   \coordinate (C5) at (8.0,4.0);
   \coordinate (D1) at (0.0,6.0);
   \coordinate (D5) at (8.0,6.0);
   \coordinate (E1) at (0.0,8.0);
   \coordinate (E2) at (2.0,8.0);
   \coordinate (E3) at (4.0,8.0);
   \coordinate (E4) at (6.0,8.0);
   \coordinate (E5) at (8.0,8.0);
   \draw (A1) -- (A5) ;
   \draw (B1) -- (B5) ;
   \draw (C1) -- (C5) ;
   \draw (D1) -- (D5) ;
   \draw (E1) -- (E5) ;
   \draw (A1) -- (E1) ;
   \draw (A2) -- (E2) ;
   \draw (A3) -- (E3) ;
   \draw (A4) -- (E4) ;
   \draw (A5) -- (E5) ;
   \draw (A4) -- (B5) ;
   \draw (A3) -- (C5) ;
   \draw (A2) -- (D5) ;
   \draw (A1) -- (E5) ;
   \draw (B1) -- (E4) ;
   \draw (C1) -- (E3) ;
   \draw (D1) -- (E2) ;
\end{tikzpicture}
\qquad\qquad
\begin{tikzpicture}[line width = 0.8pt,scale=0.5]
\coordinate (X) at (2.0,0.0);
\coordinate (Y) at (2.0,2.0);
\coordinate (Z) at (0.0,2.0);
\coordinate (S) at ($0.45*(-2,2)$);
\coordinate (T) at ($0.55*(-2,2)$);
\foreach \y in {0, 2, 4, 6}
{
\foreach \x in {0, 2, 4, 6}
{\coordinate (A) at (\x,\y);
  \coordinate (B) at ($(A)+(X)$);
  \coordinate (C) at ($(A)+(Y)$);
  \coordinate (D) at ($(A)+(Z)$);
  \coordinate (E) at ($(B)+(S)$);
  \coordinate (F) at ($(B)+(T)$);
  \draw (A) -- (B) -- (C) -- (D) -- (A) -- (C) -- (E) -- (A) -- (F) --(C);
  \draw (B) -- (E);
  \draw (D) -- (F);
  };};
\end{tikzpicture} 
\caption{Mesh $\T_h^1$ (left) and its refinement
$\T_h^2$ (right).}
 \label{fig1}
\end{figure}
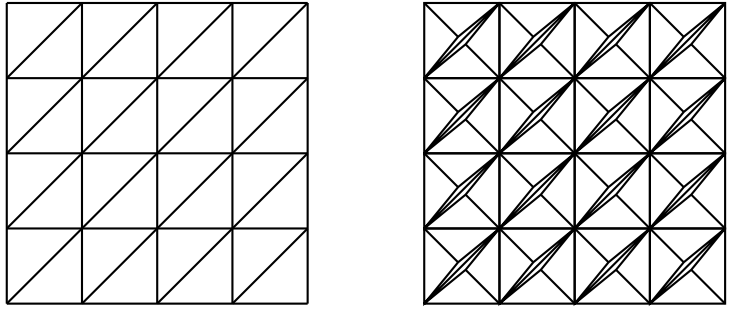

The minimum or maximum angle conditions are sufficient but not
necessary for the standard error estimation \eqref{Lag-error} or \eqref{eq2.6}.
Suppose that we have the two meshes $\T_h^1$ and $\T_h^2$ depicted in
Figure~\ref{fig1} {(where $\T_h^2$ has already been mentioned
 in \cite{HKK})}.  The mesh $\T_h^1$ consists of congruent isosceles right
triangles and satisfies the minimum angle (shape  regularity) condition.
In contrast, $\T_h^2$ is a ``refinement'' of $\T_h^1$ and contains
very flat triangles. Thus, $\T_h^2$ does not satisfy either the minimum
or maximum angle condition when the flat triangles degenerate.

As each $\tau \in \T_h^2$ is a subset of some $\tau' \in \T_h^1$,
we define the interpolation $\Pi^*$ on $\T_h^2$ as follows.
Let $\NN_h^1$ and $\NN_h^2$ be the sets of vertices of $\T_h^1$ and
$\T_h^2$, respectively.  Let $\{\varphi_i\} \subset S_h$ be the basis of
$S_h$ on $\T_h^2$ defined by
\begin{align*}
   \varphi_i(\bfx_i) = 1, \qquad
   \varphi_i(\bfx_j) = 0, \quad \bfx_i, \bfx_j \in \NN_h^2, \quad i \neq j.
\end{align*}
We define
 \begin{align*}
   w_i := \begin{cases}
     (\Pi_{\T_h^1}^1 u)(\bfx_i), &
      \text{for }\; \bfx_i \in \NN_h^2\backslash \NN_h^1, \\
     u(\bfx_i), & \bfx_i \in \NN_h^1,
   \end{cases} \qquad
  \Pi^* u := \sum_i w_i\varphi_i,
\end{align*}
where $\Pi_{\T_h^1}^1$ is the Lagrange interpolation on $\T_h^1$.
As $\T_h^1$ is a good mesh, the standard error estimation
\eqref{eq2.6} holds, and the finite element solutions converge to $u$
with $\OO(h)$.  Moreover, $\Pi^*$ inherits these properties because
$\Pi_{\T_h^1}^1u = \Pi^* u$ on $\Omega$.

This example reminds us that even if a mesh contains bad elements,
we might be able to define an interpolation $\Pi^*$ that 
satisfies \eqref{inter} if bad elements are covered by {fine}
simplices under some additional conditions. 

\section{Assumptions and the main theorem}
In this section, we introduce assumptions that ensure the standard
error estimation \eqref{inter} on meshes that contain bad elements
and explain the main result in detail.
Recall that $0 < \theta,\psi \le \pi/3$ are given fixed constants.
Let {$M$, $N$ be positive integers and $C \ge 1$}.  These
constants are also fixed throughout the paper. 

\begin{figure}[htbp]
\centering
\begin{tikzpicture}[line width = 0.5pt,scale=1]
   \coordinate (A) at (0.0,0.0);
   \coordinate (B) at (2.0,0.0);
   \coordinate (C) at (2.2,1.1);
   \coordinate (D) at (1.1,0.9);
   \coordinate (E) at (-0.2,2.0);
   \draw (A) -- (D) -- (B) -- (C) -- (D) -- (E) -- (A);
   \fill (D) circle (2.5pt);
   \coordinate [label={${T_1}$}] (F) at (0.45,0.7);
   \coordinate [label={${T_2}$}] (G) at (1.75,0.45);
\end{tikzpicture}
 \qquad\quad
\begin{tikzpicture}[line width = 0.5pt,scale=1.2]
   \coordinate (A) at (0.0,1.0);
   \coordinate (B) at (1.0,0.0);
   \coordinate (C) at (2.2,1.1);
   \coordinate (D) at (1.1,2.0);
   \draw (A) -- (B) -- (C) -- (D) -- (A);
   \fill (A) circle (2.0pt) (C) circle (2.0pt);
   \draw [line width=1.2pt] (A) -- (C);
   \coordinate [label={${T_1}$}] (F) at (1.1,0.45);
   \coordinate [label={${T_2}$}] (G) at (1.1,1.27);
\end{tikzpicture}
\qquad\quad
\begin{tikzpicture}[line width = 0.5pt,scale=1.2]
   \coordinate (A) at (0.0,1.0);
   \coordinate (B) at (1.0,0.0);
   \coordinate (C) at (2.2,1.1);
   \draw (A) -- (B) -- (C) -- (A);
   \coordinate (D) at (1.1,2.0);
   \coordinate (E) at (1.3,1.059);
   \draw (A) -- (D) -- (E);
   \fill (A) circle (2.0pt);
   \coordinate [label={${T_1}$}] (F) at (1.1,0.45);
   \coordinate [label={${T_2}$}] (G) at (0.9,1.15);
\end{tikzpicture}
\caption{Examples of $V({T_1},{T_2})$.} \label{fig2}
\end{figure}
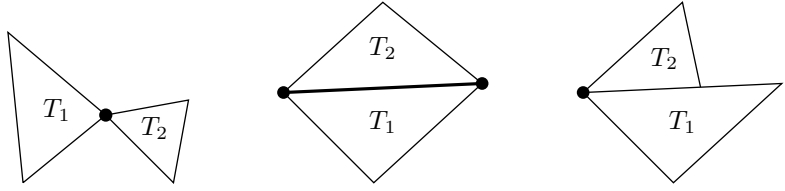

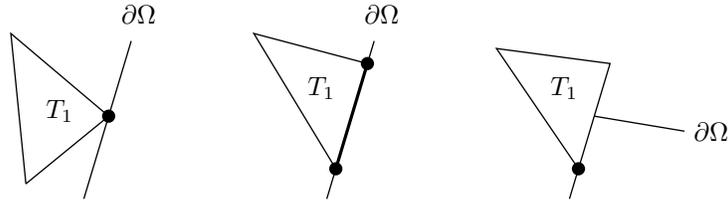
\begin{figure}[htbp]
\centering
\begin{tikzpicture}[line width = 0.5pt,scale=1]
   \coordinate (A) at (0.0,0.0);
   \coordinate (B) at (1.1,0.9);
   \coordinate (C) at (-0.2,2.0);
   \draw (A) -- (B) -- (C) -- (A);
   \fill (B) circle (2.5pt);
   \coordinate [label={${T_1}$}] (F) at (0.45,0.7);
   \coordinate (D) at (0.77,-0.2);
   \coordinate (E) at (1.4,1.9);
   \draw (D) -- (E);
   \coordinate [label={$\partial\Omega$}] (G) at (1.5,2.0);
\end{tikzpicture}
\qquad\quad
\begin{tikzpicture}[line width = 0.5pt,scale=1]
   \coordinate (A) at (0.89,0.2);
   \coordinate (B) at (1.31,1.6);
   \coordinate (C) at (-0.2,2.0);
   \coordinate (D) at (0.77,-0.2);
   \coordinate (E) at (1.4,1.9);
   \draw (D) -- (A) -- (C) -- (B) -- (E);
   \fill (B) circle (2.5pt) (A) circle (2.5pt);
   \coordinate [label={${T_1}$}] (F) at (0.7,1.0);
   \draw [line width=1.2pt] (A) -- (B);
   \coordinate [label={$\partial\Omega$}] (G) at (1.5,2.0);
\end{tikzpicture}
\qquad\quad
\begin{tikzpicture}[line width = 0.5pt,scale=1]
   \coordinate (A) at (0.89,0.2);
   \coordinate (B) at (1.31,1.6);
   \coordinate (C) at (-0.2,1.8);
   \coordinate (D) at (0.77,-0.2);
   \coordinate [label=right:{$\partial\Omega$}] (E) at (2.3,0.7);
   \coordinate (F) at (1.1,0.9);
   \draw (D) -- (A) -- (C) -- (B) -- (A);
   \draw (E) -- (F);
   \fill (A) circle (2.5pt);
   \coordinate [label={${T_1}$}] (F) at (0.7,1.0);
    \coordinate (G) at (1.5,2.0);
\end{tikzpicture}

\caption{Examples of $V_{\partial\Omega}({T_1})$.}
 \label{fig3}
\end{figure}

For a simplex $\tau$, let $v(\tau)$ be the set of vertices of $\tau$.
For $n$-simplices ${T_1}$ and ${T_2}$, we define
\begin{align*}
   V({T_1},{T_2}) := 
   \mathrm{Conv}(v({T_1})\cap v({T_2})),
   \qquad V_{\partial\Omega}({T_1}) :=
    \mathrm{Conv}(v({T_1})\cap \partial\Omega),
\end{align*}
where Conv$(S)$ denotes the convex hull of the point set $S$.
We present several examples in Figures \ref{fig2} and \ref{fig3}.
For triangles ${T_1}$ and ${T_2}$, 
$V({T_1},{T_2})$ and
$V_{\partial\Omega}({T_1})$ are drawn with thick lines and dots.

We impose the following conditions for the mesh $\T_h$.
\begin{assumption}\label{A2}
$(1)$ There exist $n$-simplices $T_1, \dots, T_m$ such that
$\theta_{T_k}^{min} \ge \psi$ and $h_{T_k} \le C h$ for
$k = 1, \dots, m$.
We assume that $\ring{T}_k \subset \Omega$, {where 
$\ring{T}_k$ is the interior set of $T_k$.} 
{The interior sets $\{\ring{T}_k\}$ may overlap, and}
the multiplicity of the overlap of $\{\ring{T}_k\}_{k=1}^m$ is
at most $M$. (Note that we do NOT assume
$T_k \in \T_h$ ($k = 1, \dots, m)$.) \\
$(2)$ The mesh $\T_h$ is expressed as a disjoint union
$\T_h = \bigcup_{k=0}^m Q_k$ such that 
{$Q_0 \subset \T_h \cap K_\theta$
and} each ``bad'' element $\tau \in \T_h \backslash K_\theta$
belongs to one of $\{Q_k\}_{k=1}^m$.
Furthermore, we assume that
\begin{align*}
   \bigcup_{\tau \in Q_k} \tau \subset \overline{T}_k, \quad
     \text{ and } \quad \#Q_k \le N \quad\text{ for }k = 1, \dots, m
     {,}
\end{align*}
{where $\#Q_k$ is the number of elements in $Q_k$.}
{
Note that $Q_0$ consists of ``good'' and ``ordinary'' elements, and
$Q_k$ $(k=1, \dots, m)$ may contain ``good'' or ``ordinary'' elements
as well as ``bad'' elements.}  \\
$(3)$ If $(v(\tau_0) \cap v(\tau_1)) \backslash v(T_k) \neq \emptyset$ for
$\tau_0 \in Q_0$, $\tau_1 \in Q_k$, $1 \le k \le m$,
then we assume that $\theta_{\tau_0}^{min} \ge {\theta}$. 
(Note that if $n = 3$, this condition is
automatically satisfied {because all elements belong to
$Q_0$ are ``good''}.)  \\
$(4)$ If $v(\tau_1) \cap v(\tau_2) \neq \emptyset$ for
$\tau_1 \in Q_k$, $\tau_2 \in Q_\ell$, $1 \le k < \ell \le m$,
then we assume $v(\tau_1) \cap v(\tau_2) \subset V(T_k,T_\ell)$. \\
$(5)$ If $v(\tau_1) \cap \partial\Omega \neq \emptyset$ for
$\tau_1 \in Q_k$, $1 \le k \le m$, then
we assume $v(\tau_1) \cap \partial\Omega \subset V_{\partial\Omega}(T_k)$.
\end{assumption}

Let $\{\varphi_i\} \subset S_h$ be the basis of $S_h$ on the mesh $\T_h$.
\begin{definition}
For $u \in H^2(\Omega)$, we define the interpolation $\Pi^* u$ on $\T_h$ by
\begin{align*}
   w_i & := \begin{cases}
     (\Pi_{T_k}^1 u)(\bfx_i), &
      \text{if $\bfx_i \in v(\tau)$ for some $\tau \in Q_k$} \\
     u(\bfx_i), & \text{otherwise}
   \end{cases} \\
  \Pi^* u & := \sum_i w_i\varphi_i,
\end{align*}
where {$\bfx_i$ are the vertices of the mesh  $\T_h$ and }
$\Pi_{T_k}^1$ is the Lagrange interpolation on $T_k$.
\end{definition}

\textit{Remark 5.} Assumption~\ref{A2} (1), (2) requires that
{each bad element in $\T_h$ belongs to one of the disjoint 
element clusters $Q_k$ ($k = 1, \dots, m$) and each $Q_k$ is
virtually covered by a simplex $T_k$ ($k = 1, \dots, m$)}.
The interpolation $\Pi^* u$ is defined by $\Pi_{T_k}^1 u$ on each
$T_k$. 

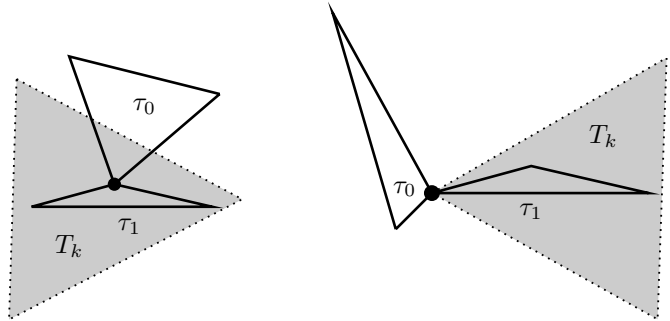
\begin{figure}[htbp]
\centering
\begin{tikzpicture}[line width = 1.0pt,scale=1.0]
  % \draw [help lines] (-2,-2) grid (3,3);
   \coordinate (X) at (-0.2,1.7);
   \coordinate (Y) at (2.8,0.1);
   \coordinate (Z) at (-0.3,-1.5);
   \fill [black!20] (X) -- (Y) -- (Z) -- (X) -- cycle;
   \coordinate (A) at (2.5,1.5);
   \coordinate (B) at (2.4,0.0);
   \coordinate (C) at (1.1,0.3);
   \coordinate (D) at (0.0,0.0);
   \coordinate (E) at (0.5,2.0);
   \draw (A) -- (C) -- (E) -- (A);
   \draw (D) -- (B) -- (C) -- (D);
   \fill (C) circle (2.5pt);
   \coordinate [label={$\tau_0$}] (F) at (1.5,1.1);
   \coordinate [label={$\tau_1$}] (G) at (1.3,-0.5);
    \draw [dotted, line width = 0.7pt](X) -- (Y) -- (Z) -- (X);
   \coordinate [label={$T_k$}] (G) at (0.5,-0.8); 
\end{tikzpicture}
 \qquad\quad
\begin{tikzpicture}[line width = 1.0pt,scale=1.2]
   \coordinate (X) at (0.0,0.0);
   \coordinate (Y) at (2.6,1.5);
   \coordinate (Z) at (2.5,-1.4);
   \fill [black!20] (X) -- (Y) -- (Z) -- (X) -- cycle;
   \coordinate (A) at (-0.4,-0.4);
   \coordinate (B) at (2.4,0.0);
   \coordinate (C) at (1.1,0.3);
   \coordinate (D) at (0.0,0.0);
   \coordinate (E) at (-1.1,2.0);
   \draw (A) -- (D) -- (E) -- (A);
   \draw (D) -- (B) -- (C) -- (D);
   \fill (D) circle (2.5pt);
   \coordinate [label={$\tau_0$}] (F) at (-0.3,-0.15);
   \coordinate [label={$\tau_1$}] (G) at (1.1,-0.4);
    \draw [dotted, line width = 0.7pt](X) -- (Y) -- (Z) -- (X);
   \coordinate [label={$T_k$}] (G) at (1.9,0.4); 
\end{tikzpicture}
\caption{{Two situations regarding Assumption~3 (3).}} \label{fig4}
\end{figure}

{
To explain the meaning of Assumption~\ref{A2} (3), we consider the
situation that $\tau_0 \in Q_0$ and $\tau_1 \in Q_k$ as is
depicted in Figure~\ref{fig4}. 
Then, $\tau_1$ is virtually covered by $T_k$ (with other elements in
$Q_k$, if $\{\tau_1\} \subsetneqq Q_k$).  Let
$\bfx \in v(\tau_0) \cap v(\tau_1) \backslash v(T_k)$.
Figure~\ref{fig4} (left) depicts such a situation.
As the interpolation $\Pi^*u$ is defined by $\Pi_{T_k}^1 u$,
we have
$(\Pi^* u)(\bfx) = (\Pi_{T_k}^1 u)(\bfx) \neq (\Pi_{\tau_0}^1u)(\bfx)$
in general in this case.  To absorb the effect of
$(\Pi^* u)(\bfx) \neq (\Pi_{\tau_0}^1u)(\bfx)$,
we use the shape regularity condition \eqref{SR} for $\tau_0$ and
the Sobolev embedding inequality \eqref{eq2.8} for
$T_k$.   Therefore,
we require that $\tau_0$ satisfies the minimum angle
condition. 

If $\bfx \in v(T_k)$, in contrast, we have 
$(\Pi^* u)(\bfx) = (\Pi_{T_k}^1 u)(\bfx) = (\Pi_{\tau_0}^1u)(\bfx)$,
and we do not need to impose such a condition
(see Figure~\ref{fig4} (right)).}

\begin{figure}[htbp]
\centering
\begin{tikzpicture}[line width = 1.0pt,scale=1.0]
   \coordinate (X) at (-0.2,1.7);
   \coordinate (Y) at (2.8,0.1);
   \coordinate (Z) at (-0.3,-1.5);
   \coordinate (S) at (0.4,1.2);
   \coordinate (T) at (-2.8,0.1);
   \coordinate (U) at (-0.2,-1.3);
   \fill [black!20] (X) -- (Y) -- (Z) -- (X) -- cycle;
   \fill [black!20] (S) -- (T) -- (U) -- (S) -- cycle;
   \coordinate (A) at (-2.1,-0.15);
   \coordinate (B) at (2.4,0.0);
   \coordinate (C) at (1.1,0.4);
   \coordinate (D) at (0.0,0.0);
   \coordinate (E) at (-0.9,0.4);
   \draw (A) -- (D) -- (E) -- (A);
   \draw (D) -- (B) -- (C) -- (D);
   \fill (D) circle (2.5pt);
   \coordinate [label={$\tau_1$}] (F) at (-0.9,-0.1);
   \coordinate [label={$\tau_2$}] (G) at (1.2,-0.05);
   \draw [dotted, line width = 0.7pt](X) -- (Y) -- (Z) -- (X);
   \draw [dotted, line width = 0.7pt](S) -- (T) -- (U) -- (S);
   \coordinate [label={$T_\ell$}] (G) at (0.9,0.4); 
   \coordinate [label={$T_k$}] (G) at (-0.9,-0.75); 
\end{tikzpicture}
\caption{{A situation that makes $\Pi^* u$ be not well-defined.}}
\label{fig5}
\end{figure}
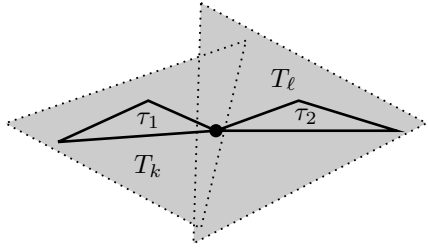

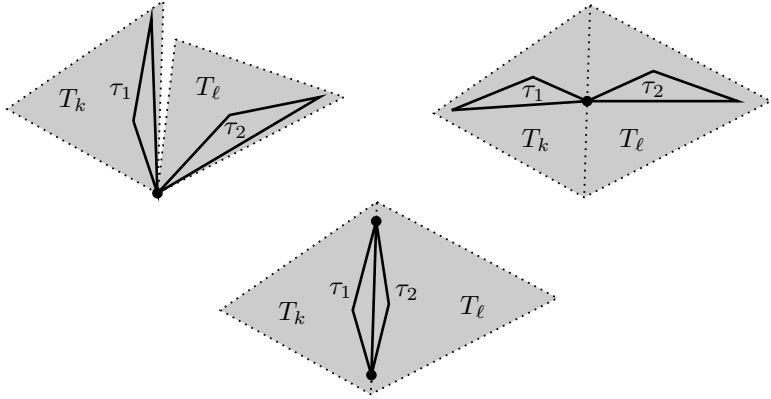
\begin{figure}[htbp]
\centering
\begin{tikzpicture}[line width = 1.0pt,scale=0.8]
   \coordinate (X) at (0.1,3.2);
   \coordinate (Y) at (3.1,1.6);
   \coordinate (Z) at (0.0, 0.0);
   \coordinate (U) at (-2.5,1.4);
   \coordinate (S) at (0.3,2.56);
   \fill [black!20] (X) -- (Z) -- (U) -- (X) -- cycle;
   \fill [black!20] (S) -- (Y) -- (Z) -- (S) -- cycle;
   \coordinate (A) at (-0.4,1.2);
   \coordinate (B) at (2.7,1.6);
   \coordinate (C) at (1.2,1.3);
   \coordinate (E) at (-0.1,2.9);
   \draw (A) -- (Z) -- (E) -- (A);
   \draw (Z) -- (B) -- (C) -- (Z);
   \fill (Z) circle (2.5pt);
   \coordinate [label={$\tau_1$}] (F) at (-0.6,1.4);
   \coordinate [label={$\tau_2$}] (G) at (1.3,0.75);
   \draw [dotted, line width = 0.7pt] (S) -- (Y) -- (Z) -- (S);
   \draw [dotted, line width = 0.7pt] (Z) -- (U) -- (X) -- (Z);
   \coordinate [label={$T_k$}] (H) at (-1.4,1.2); 
   \coordinate [label={$T_\ell$}] (I) at (0.85,1.4); 
\end{tikzpicture}
 \qquad\quad
\begin{tikzpicture}[line width = 1.0pt,scale=0.8]
   \coordinate (X) at (0.1,3.2);
   \coordinate (Y) at (3.1,1.6);
   \coordinate (Z) at (0.0, 0.0);
   \coordinate (U) at (-2.5,1.4);
   \fill [black!20] (X) -- (Y) -- (Z) -- (U) -- (X) -- cycle;
   \coordinate (A) at (-2.2,1.45);
   \coordinate (B) at (2.55,1.6);
   \coordinate (C) at (1.15,2.1);
   \coordinate (D) at (0.05,1.6);
   \coordinate (E) at (-0.85,2.0);
   \draw (A) -- (D) -- (E) -- (A);
   \draw (D) -- (B) -- (C) -- (D);
   \fill (D) circle (2.5pt);
   \coordinate [label={$\tau_1$}] (F) at (-0.84,1.43);
   \coordinate [label={$\tau_2$}] (G) at (1.15,1.5);
   \draw [dotted, line width = 0.7pt] (X) -- (Y) -- (Z) -- (U) -- (X) -- (Z);
   \coordinate [label={$T_k$}] (H) at (-0.8,0.6); 
   \coordinate [label={$T_\ell$}] (I) at (0.8,0.6); 
\end{tikzpicture} \\
\begin{tikzpicture}[line width = 1.0pt,scale=0.8]
   \coordinate (X) at (0.1,3.2);
   \coordinate (Y) at (3.1,1.6);
   \coordinate (Z) at (0.0, 0.0);
   \coordinate (U) at (-2.5,1.4);
   \fill [black!20] (X) -- (Y) -- (Z) -- (U) -- (X) -- cycle;
   \draw [dotted, line width = 0.7pt] (X) -- (Y) -- (Z) -- (U) -- (X) -- (Z);
   \coordinate (A) at (0.01,0.32);
   \coordinate (B) at (0.09,2.88);
   \coordinate (C) at (0.3,1.5);
   \coordinate (D) at (-0.3,1.4);
   \fill (A) circle (2.5pt);
   \fill (B) circle (2.5pt);
   \draw (A) -- (C) -- (B) -- (D) -- (A) -- (B);
   \coordinate [label={$\tau_1$}] (F) at (-0.5,1.43);
   \coordinate [label={$\tau_2$}] (G) at (0.6,1.4);
   \coordinate [label={$T_k$}] (H) at (-1.3,1.0); 
   \coordinate [label={$T_\ell$}] (I) at (1.7,1.1); 
\end{tikzpicture}
\caption{{Three situations that satisfy Assumption~\ref{A2}(4).}} \label{fig6}
\end{figure}

Suppose now that 
$v(\tau_1) \cap v(\tau_2) \neq \emptyset$,
$\tau_1 \in Q_k$, $\tau_2 \in Q_\ell$, 
($1 \le k < \ell \le m$) and
$\bfx \in v(\tau_1) \cap v(\tau_2)$ is an
inner point of $\Omega$.  
If $(\Pi_{T_k}^1 u)(\bfx) \neq (\Pi_{T_\ell}^1u) (\bfx)$, then
$\Pi^* u$ is not well-defined.  This case may occur if $\bfx$, $T_k$,
and $T_\ell$ are in the situation depicted in Figure~\ref{fig5}.
Assumption~\ref{A2} (4) is imposed to avoid such a situation.
Three situations that satisfy Assumption~\ref{A2} (4) are depicted
in Figure~\ref{fig6}.

\begin{figure}[htbp]
\centering
\begin{tikzpicture}[line width = 1.0pt,scale=0.8]
   \coordinate (A) at (0.2,1.6);
   \coordinate (Y) at (-2.0, 3.2);
   \coordinate (Z) at (-2.5,0.0);
   \fill [black!20] (A) -- (Y) -- (Z) -- (A) -- cycle;
   \coordinate (S) at (-0.04, -0.32);
   \coordinate (T) at (0.4, 3.2);
   \draw [dotted, line width = 0.7pt] (A) -- (Y) -- (Z) -- (A);
   \draw [line width = 0.7pt] (S) -- (T);
    \coordinate (B) at (-0.9,1.9);
   \coordinate (C) at (-2.1,1.5);
   \draw (A) -- (B) -- (C) -- (A);
   \fill (A) circle (2.5pt);
   \coordinate [label={$\tau_1$}] (D) at (-1.0,1.0);
   \coordinate [label={$T_k$}] (E) at (-1.6,2.0); 
   \coordinate [label={$\partial\Omega$}] (F) at (-0.4,-0.2); 
 \end{tikzpicture}
\qquad \quad
\begin{tikzpicture}[line width = 1.0pt,scale=0.8]
   \coordinate (X) at (0.4,3.2);
   \coordinate (Y) at (0.0,0.0);
   \coordinate (Z) at (-2.5,1.4);
   \fill [black!20] (Y) -- (X) -- (Z) -- (Y) -- cycle;
   \coordinate (S) at (-0.05, -0.4);
   \coordinate (T) at (0.45, 3.6);
   \draw [dotted, line width = 0.7pt] (Y) -- (Z) -- (X);
   \draw [line width = 0.7pt] (S) -- (T);
   \coordinate (A) at (0.2,1.6);
   \coordinate (B) at (-0.1,2.7);
   \coordinate (C) at (-0.4,0.33);
   \draw (A) -- (B) -- (C) -- (A);
   \fill (A) circle (2.5pt);
   \coordinate [label={$\tau_1$}] (D) at (-0.6,1.4);
   \coordinate [label={$T_k$}] (E) at (-1.4,1.2); 
   \coordinate [label={$\partial\Omega$}] (F) at (0.4,-0.2); 
 \end{tikzpicture}
 \qquad\quad
\begin{tikzpicture}[line width = 1.0pt,scale=0.8]
   \coordinate (X) at (1.5,3.2);
   \coordinate (Y) at (3.6,0.6);
   \coordinate (Z) at (0.0, 0.0);
   \coordinate (S) at (-0.2,-0.03333);
   \coordinate (T) at (3.8, 0.0);
   \fill [black!20] (X) -- (Y) -- (Z) -- (X) -- cycle;
   \draw [dotted, line width = 0.7pt] (Y) -- (X) -- (Z);
   \draw [line width = 0.7pt] (S) -- (Y) -- (T);
   \coordinate (A) at (1.8,0.3);
   \coordinate (B) at (0.3,0.5);
   \draw (Y) -- (A) -- (B) -- (Y);
   \fill (Y) circle (2.5pt);
   \fill (A) circle (2.5pt);
   \coordinate [label={$\tau_1$}] (F) at (1.7,0.6);
   \coordinate [label={$T_k$}] (H) at (1.7,1.8); 
   \coordinate [label={$\partial\Omega$}] (I) at (4.1,0.1); 
\end{tikzpicture}
\caption{{Three situations that satisfy Assumption~\ref{A2}(5).}} \label{fig7}
\end{figure}
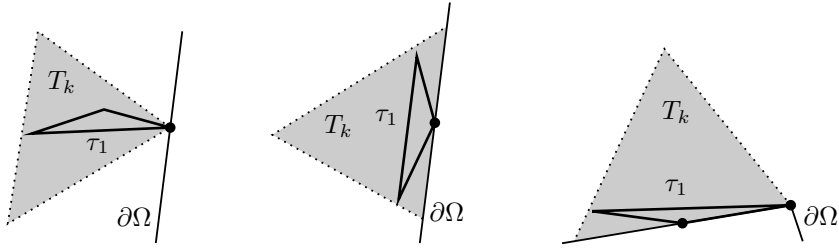

{Suppose that $\tau_1 \in Q_k$ ($1 \le k \le m$) with
$v(\tau_1) \cap \partial\Omega \neq \emptyset$. The interpolation
$\Pi^*u$ should then satisfy the boundary condition 
$\Pi^*u = 0$ on $\partial\Omega$.  Assumption~\ref{A2}(5) is imposed
to ensure that this condition is fulfilled.
Three situations that satisfy Assumption~\ref{A2}(5) are depicted
in Figure~\ref{fig7}. }

\setcounter{theorem}{5}
\noindent
\begin{theorem}\label{thm5}
Suppose that the mesh $\T_h$ satisfies Assumption~\ref{A2}.
The interpolation $\Pi^* u$ is then well defined for
$u \in H^2(\Omega)$, and we have
\begin{align}
  & |u - \Pi^* u|_{H^1(\Omega)} \le  E h |u|_{H^2(\Omega)},
  \qquad  \forall u \in H^2(\Omega), \label{eq3.1} \\
  & E := \left(A_{\theta,\psi}^2 (n+2 + C^2M) + 
    \frac{2(n+1)(n+2) \pi B_\psi^2 C^{4-n} D_{{\theta}} MN}{n{\theta}}\right)^{1/2}.
  \label{eq3.2}
\end{align}
\end{theorem}
\proof Define $W \subset Q_0$ by
$W := \{\tau \in Q_0 \mid \theta_\tau^{min} \ge {\theta}\}$.
We decompose $|u - \Pi^*u|_{H^1(\Omega)}^{{2}}$ as
\begin{align*}
  |u - \Pi^* u|_{H^1(\Omega)}^2
    = \sum_{\tau \in W}
    |u - \Pi^* u|_{H^1(\tau)}^2
    & + \sum_{\tau \in Q_0\backslash W}
    |u - \Pi^* u|_{H^1(\tau)}^2 \\
    & + \sum_{k=1}^m\sum_{\tau \in Q_k}
    |u - \Pi^* u|_{H^1(\tau)}^2.
\end{align*}
We estimate each term on the right-hand side.  As
\begin{align*}
  \Pi^* u = \Pi_\tau^1 u - (\Pi_\tau^1 u - \Pi^* u)
          = \Pi_\tau^1 u - \sum_{\bfx_i \in v(\tau)} 
            (u(\bfx_i) - w_i) \varphi_i 
\end{align*}
for $\tau \in W$, we have
\begin{align*}
  u - \Pi^* u  = u - \Pi_\tau^1 u + \sum_{\bfx_i \in v(\tau)} 
            (u(\bfx_i) - w_i) \varphi_i.
\end{align*}
As $u(\bfx_i) \neq w_i$ occurs to at most $n+1$ vertices
of $\tau$,
by the triangle and Cauchy--Schwarz inequalities, we have
\begin{align*}
   \left|u - \Pi^* u\right|_{H^1(\tau)}^2 & \le
    (n+2)\left|u - \Pi_\tau^1 u\right|_{H^1(\tau)}^2
   + (n+2)\sum_{\bfx_i \in v(\tau)}
    \left|u(\bfx_i) - w_i\right|^2 |\varphi_i|_{H^1(\tau)}^2.
\end{align*}
Hence, we have
\begin{align*}
 \sum_{\tau \in W} \left|u - \Pi^* u\right|_{H^1(\tau)}^2
 & \le (n+2) \sum_{\tau \in W} \left|u - \Pi_\tau u\right|_{H^1(\tau)}^2 \\
 & \qquad + (n+2) \sum_{\tau \in W} \sum_{\bfx_i \in v(\tau)}
    \left|u(\bfx_i) - w_i\right|^2 |\varphi_i|_{H^1(\tau)}^2 \\
 (\text{by \eqref{eq2.6}}) \quad
 & \le (n+2)A_{\theta,\psi}^2 h^2 \sum_{\tau \in W}|u|_{H^2(\tau)}^2 \\
 & \qquad + (n+2)\sum_{\tau \in W} \sum_{\bfx_i \in v(\tau)}
    \left|u(\bfx_i) - w_i\right|^2 |\varphi_i|_{H^1(\tau)}^2 \\
 & \le (n+2)A_{\theta,\psi}^2 h^2 \sum_{\tau \in W}|u|_{H^2(\tau)}^2 \\
 & \qquad + (n+2)\sum_{\tau \in W} \frac{D_{{\theta}} h_\tau^{n-2}}{n(n-1)}
     \sum_{\bfx_i \in v(\tau)}\left|u(\bfx_i) - w_i\right|^2.
\end{align*}
Here, we used the fact that
\begin{align}
  |\nabla \varphi_i| & \le \frac{1}{e_\tau}
   = \frac{|{L}_\tau|}{n|\tau|}
   \le {\frac{h_\tau^{n-1}}{n(n-1)|\tau|}}, \notag \\
  |\varphi_i|_{H^1(\tau)}^2 & \le \frac{|\tau|}{e_\tau^2}
    \le {  \frac{h_\tau^{n-1}}{n(n-1)e_\tau}      }
  \le \frac{h_\tau^{n-1}}{n(n-1)\rho_\tau}
    \le \frac{D_{{\theta}} h_\tau^{n-2}}{n(n-1)}
  \quad (\text{by \eqref{SR}}),
  \label{assump3-4-1}
\end{align}
where ${L}_\tau$ is the longest edge or widest facet
of $\tau$, and $e_\tau$ is the height of $\tau$ with respect
to ${L}_\tau$. 

For the last term in the above inequality, we note that
{\allowdisplaybreaks
\begin{align*}
  (n+2) \sum_{\tau \in W} & \frac{D_{{\theta}} h_\tau^{n-2}}{n(n-1)}
  \sum_{\bfx_i \in v(\tau)} |u(\bfx_i) - w_i|^2 \\
  & \le (n+2) \frac{D_{{\theta}} h^{n-2}}{n(n-1)} 
    \sum_{\tau \in W}\sum_{k=1}^m \sum_{\beta \in Q_k} 
   \sum_{\bfx_i \in v(\tau) \cap v(\beta)} |u(\bfx_i) - w_i|^2 \\
  & \le (n+2) \frac{D_{{\theta}} h^{n-2}}{n(n-1)} \sum_{\tau \in W}
   \sum_{k=1}^m \sum_{\beta \in Q_k} 
   \sum_{\bfx_i \in v(\tau) \cap v(\beta)} 
     \|u - \Pi_{T_k}^1 u\|_{L^\infty(T_k)}^{{2}} \\
 \text{(by \eqref{eq2.8})} \;
  & \le (n+2) \frac{D_{{\theta}} h^{n-2}}{n(n-1)} \sum_{\tau \in W}
   \sum_{k=1}^m \sum_{\beta \in Q_k} 
   \sum_{\bfx_i \in v(\tau) \cap v(\beta)} 
     B_\psi^2 h_{T_k}^{4-n}|u|_{H^2(T_k)}^2 \\
  \text{(by Assump{tion~}\ref{A2}(1))} \;
  & \le (n+2) \frac{D_{{\theta}} h^{n-2}}{n(n-1)} \sum_{\tau \in W}
   \sum_{k=1}^m \sum_{\beta \in Q_k} 
   \sum_{\bfx_i \in v(\tau) \cap v(\beta)} 
     B_\psi^2 C^{4-n}h^{4-n}|u|_{H^2(T_k)}^2 \\
  & = (n+2) \frac{B_\psi^2 C^{4-n} D_{{\theta}} h^{2}}{n(n-1)}
   \sum_{k=1}^m |u|_{H^2(T_k)}^2 \sum_{\beta \in Q_k} 
    \sum_{\tau \in W}
   \sum_{\bfx_i \in v(\tau) \cap v(\beta)} 1.
\end{align*}
}

Note that, for $\beta \in Q_k$ and each vertex $\bfx_i \in v(\beta)$,
there are at most
\begin{align}
   \left\lfloor \frac{2(n-1)\pi}{{\theta}} \right\rfloor
   \label{assump3-4-2}
\end{align}
simplices such that $\tau \in W$ and
$\bfx_i \in v(\tau) \cap v(\beta)$ 
{because of Assumption~\ref{A2}(4)}. 
Here, $\lfloor x\rfloor$ is the floor function for $x \in \R$.
Therefore, we have
\begin{align*}
   (n+2) &\frac{B_\psi^2 C^{4-n} D_{{\theta}} h^{2}}{n(n-1)}
    \sum_{k=1}^m |u|_{H^2(T_k)}^2 \sum_{\beta \in Q_k} 
    \sum_{\tau \in W}
   \sum_{\bfx_i \in v(\tau) \cap v(\beta)} 1  \\
  & \le (n+2) \frac{B_\psi^2 C^{4-n} D_{{\theta}} h^{2}}{n(n-1)}
   \sum_{k=1}^m |u|_{H^2(T_k)}^2 \sum_{\beta \in Q_k} 
    (n+1) \left\lfloor \frac{2(n-1)\pi}{{\theta}} \right\rfloor \\
  \text{(by Assump{tion~}\ref{A2}(2))} \quad
   & \le \frac{2(n+1)(n+2) \pi B_\psi^2 C^{4-n} D_{{\theta}} N h^{2}}{n {\theta}}
   \sum_{k=1}^m |u|_{H^2(T_k)}^2 \\
   \text{(by Assump{tion~}\ref{A2}(1))} \quad
   & \le \frac{2(n+1)(n+2) \pi B_\psi^2 C^{4-n} D_{{\theta}} M N h^{2}}{n {\theta}}
   |u|_{H^2(\Omega)}^2.
\end{align*}
Combining these inequalities, we obtain
\begin{align*}
  \sum_{\tau \in W} \left|u - \Pi^* u\right|_{H^1(\tau)}^2
   & \le (n+2)A_{\theta,\psi}^2 h^2 \sum_{\tau \in W}|u|_{H^2(\tau)}^2 \\
  & \quad + \frac{2(n+1)(n+2) \pi B_\psi^2 C^{4-n} D_{{\theta}} M N h^{2}}{n {\theta}}
   |u|_{H^2(\Omega)}^2.
\end{align*}

Suppose now that $\tau \in Q_0 \backslash W$.  By
Assumption~\ref{A2}(3), we have $\Pi^* u = \Pi_\tau^1 u$ on
$\tau$ as is pointed out in Remark~5.  Hence, we have
\begin{align*}
   \sum_{\tau \in Q_0\backslash W}
    |u - \Pi^* u|_{H^1(\tau)}^2
   = \sum_{\tau \in Q_0\backslash W}
    |u - \Pi_\tau^1 u|_{H^1(\tau)}^2 
    \le A_{\theta,\psi}^2 h^2 \sum_{\tau \in Q_0\backslash W}
     |u|_{H^2(\tau)}^2
\end{align*}
by \eqref{eq2.6}. 

As $\Pi^* u = \Pi_{T_k}^1u$ on $\tau \in Q_k$, $k = 1, \dots, m$
by definition, we have
\begin{align*}
   \sum_{k=1}^m\sum_{\tau \in Q_k}
    |u - \Pi^* u|_{H^1(\tau)}^2 & \le
   \sum_{k=1}^m  |u - \Pi_{T_k}^1 u|_{H^1(T_k)}^2 \\
  \text{(by \eqref{eq2.6})} \quad 
 & \le  \sum_{k=1}^m  A_{\theta,\psi}^2 h_{T_k}^2 
     |u|_{H^2(T_k)}^2 \\
  \text{(by Assump{tion~}\ref{A2}(1))} \quad
 & \le   A_{\theta,\psi}^2 C^2 M h^2 
     |u|_{H^2(\Omega)}^2.
\end{align*}

Gathering these three inequalities, we obtain
\begin{align*}
    |u - \Pi^* u|_{H^1(\Omega)}^2
    & = \sum_{\tau \in W}
    |u - \Pi^* u|_{H^1(\tau)}^2
     + \sum_{\tau \in Q_0\backslash W}
    |u - \Pi^* u|_{H^1(\tau)}^2 \\
   & \hspace{3.4cm}
     + \sum_{k=1}^m\sum_{\tau \in Q_k}
    |u - \Pi^* u|_{H^1(\tau)}^2 \\
  & \le  (n+2)A_{\theta,\psi}^2 h^2 \sum_{\tau \in W}|u|_{H^2(\tau)}^2 \\
  & \quad + \frac{2(n+1)(n+2) \pi B_\psi^2 C^{4-n}
 D_{{\theta}} M N h^{2}}{n {\theta}}
   |u|_{H^2(\Omega)}^2 \\
  & \quad  + A_{\theta,\psi}^2 h^2 \sum_{\tau \in Q_0\backslash W}
     |u|_{H^2(\tau)}^2
   + A_{\theta,\psi}^2 C^2 M h^2 
     |u|_{H^2(\Omega)}^2 \\
  & \hspace{-10mm} \le \left(
A_{\theta,\psi}^2 (n+2 + C^2M)
  + \frac{2(n+1)(n+2) \pi B_\psi^2 C^{4-n} D_{{\theta}} M N}{n {\theta}}
  \right) h^2|u|_{H^2(\Omega)},
\end{align*}
and Theorem~\ref{thm5} is proved. 
\endproof

\section{Simpler formulation for the case $n=2$.}
{For the case $n = 2$, we introduce an assumption simpler
than Assumption~\ref{A2} that is still practical
and applicable to many meshes.} 

\begin{assumption}\label{A5}
$(1)$ There exist similar isosceles triangles $T_1, \dots, T_m$
whose inner angles are ${\phi}$, ${\phi}$, and
 $\pi - 2{\phi}$.
We suppose that $\ring{T}_k \subset \Omega$, $k = 1, \dots, m$, 
$\ring{T}_i \cap \ring{T}_j = \emptyset$ ($i \neq j$,
$i,j = 1, \dots, m$).  \\
$(2)$ If $\theta_\tau^{max} > \pi - {\phi}$ for $\tau \in \T_h$,
then $\ring{\tau} \subset \ring{T_k}$ with some $k$, and
the longest edges of $\tau$ and $T_k$ coincide.
Let $\bfx \in v(\tau) \backslash v(T_k)$.
All $\tau \in \T_h$ that share $\bfx$ must satisfy
$\theta_{\tau}^{min} \ge {\phi}/2$.
\end{assumption}
\begin{center}
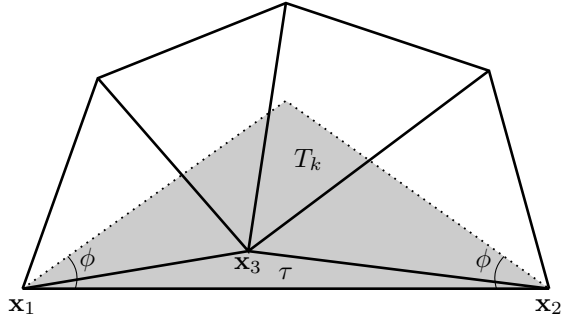
\begin{figure}[htbp]
\begin{tikzpicture}[line width = 1pt,scale=1.0]
   \coordinate [label=below:{$\bfx_1$}](A) at (0.0,0.0);
   \coordinate [label=below:{$\bfx_2$}](B) at (7.0,0.0);
   \coordinate (C) at (3.0,0.5);
   \coordinate (D) at (1.0,2.8);
   \coordinate (E) at (3.5,3.8);
   \coordinate (F) at (6.2,2.9);
   \coordinate (G) at (3.5,2.5);
   \fill [black!20] (A)--(B)--(G)--(A)--cycle;
   \draw (A) -- node[above]{$\tau$}(B) ;
   \draw (B) -- (C) ;
   \draw (C) -- (A) ;
   \draw (A) -- (D) -- (E) -- (F) -- (B);
   \draw (D) -- (C) -- (E);
   \draw (C) -- (F);
   \draw [dotted, line width = 0.7pt] (A) -- (G) -- (B);
   \coordinate (H) at (0.7,0.0);
   \coordinate (I) at (0.6,0.429);
  \draw [bend right,thin] (H) to
      node[pos=0.9,right]{${\phi}$} (I) ;
   \coordinate (J) at (6.3,0.0);
   \coordinate (K) at (6.4,0.429);
  \draw [bend left,thin] (J) to
    node[pos=0.85,left]{${\phi}$} (K) ;
   \coordinate [label=below:{$\bfx_3$}](H) at (3.0,0.55);
   \coordinate [label=below:{$T_k$}](I) at (3.8,2.0);
\end{tikzpicture} 
\caption{A bad element $\tau = \triangle \bfx_1\bfx_2\bfx_3$
  is covered virtually by a triangle $T_k$. They share the edge
 $\bfx_1\bfx_2$.  The minimum inner angles of all elements except 
 $\tau$ that share $\bfx_3$ are greater than or equal to ${\phi}/2$.}
\label{fig20}
\end{figure}
\end{center}
Figure~\ref{fig20} depicts a typical situation.
{Assumption~\ref{A5} is a special case of Assumption~\ref{A2}.
 In fact,
if} Assumption~\ref{A5} is satisfied, Assumption~\ref{A2} is satisfied
with $M = N = 1$, $C = 1$, {$\theta = \phi/2$}, and
$\psi = {\phi}$.
Hence, we have the following corollary.
\begin{corollary}\label{cor7}
 Suppose that the mesh $\T_h$ satisfies Assumption~\ref{A5}.
Then, the interpolation $\Pi^* u$ is well defined for
$u \in H^2(\Omega)$, and we have 
\begin{align*}
  |u - \Pi^* u|_{H^1(\Omega)} \le  E h |u|_{H^2(\Omega)},
  \qquad  \forall u \in H^2(\Omega),
\end{align*}
where $E$ is a constant that depends only on ${\phi}$ and
is independent of $h > 0$.
\end{corollary}

\begin{figure}[htbp]
\begin{tikzpicture}[line width = 1pt,scale=0.8]
   \coordinate (A) at (0.0,0.0);
   \coordinate (B) at (1.0,0.2);
   \coordinate (C) at (2.0,-0.2);
   \coordinate (D) at (3.0,0.0);
   \coordinate (E) at (4.0,0.2);
   \coordinate (F) at (5.0,-0.2);
   \coordinate (G) at (6.0,0.0);
   \coordinate (H) at (7.0,0.2);
   \coordinate (I) at (8.0,-0.2);
   \coordinate (J) at (9.0,0.0);
   \fill [black!20] (A)--(3.0,1.5)-- (G) -- (3.0,-1.5)--(A)--cycle;
   \fill [black!20] (G) -- (9.0,-1.5)--(9.0,1.5)--cycle;
   \draw [dotted, line width = 0.7pt] (3.0,-1.5) -- (A) --(3.0,1.5)--(9.0,-1.5) -- (9.0,1.5) -- (3.0,-1.5) -- (3.0,1.5);
   \draw (A) -- (B) -- (F) -- (H) -- (J) -- (I) -- (E) -- (C) -- (A);
   \draw (B) -- (C);
   \draw (E) -- (F);
   \draw (H) -- (I);
\end{tikzpicture} 
\caption{{Bad elements that satisfy Assumption~\ref{A2} but
do not satisfy Assumption~\ref{A5}. }}
\label{fig21}
\end{figure}
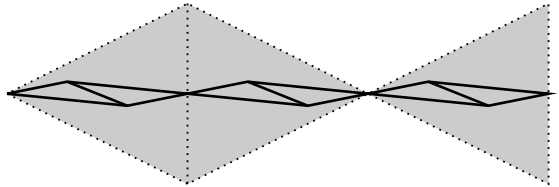

Figure~\ref{fig21} shows a mesh that satisfies Assumption~\ref{A2} but
does not satisfy Assumption~\ref{A5}.

\section{Comparison of the results with the results of prior works}
In this section, we compare our approach with approaches adopted in
prior works.  As stated, the approach adopted in this paper was
initiated by Kuc\v{e}ra \cite{Kucera} 
and Duprez, Lleras, and Lozinski \cite{DLL}.

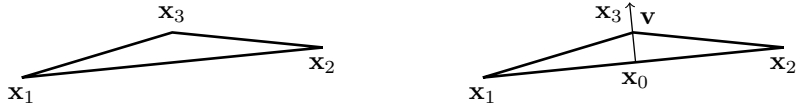
\begin{figure}[htbp]
\centering
\begin{tikzpicture}[line width = 1.0pt,scale=1.0]
   \coordinate [label=below:$\bfx_1$] (A) at (0.0,0.0);
   \coordinate [label=below:$\bfx_2$] (B) at (4.0, 0.4);
   \coordinate [label=above:$\bfx_3$] (C) at (2.0,0.6);
   \draw  (A) -- (B) -- (C) -- (A);
 \end{tikzpicture}
\qquad \qquad
\begin{tikzpicture}[line width = 1.0pt,scale=1.0]
   \coordinate [label=below:$\bfx_1$] (A) at (0.0,0.0);
   \coordinate [label=below:$\bfx_2$] (B) at (4.0, 0.4);
   \coordinate [label=above left:$\bfx_3$] (C) at (2.0,0.6);
   \draw  (A) -- (B) -- (C) -- (A);
   \coordinate [label=below:$\bfx_0$] (D) at (2.03, 0.203);
   \coordinate [label=below right:$\bfv$] (E) at (1.95, 1.0);
   \draw [->, line width = 0.5] (D) -- (E);
 \end{tikzpicture}
\caption{A bad element and the modified Lagrange
 interpolation.}
 \label{fig8}
\end{figure}

First, we explain the approach of Kuc\v{e}ra \cite{Kucera}.
Suppose that there exists a bad triangle $\tau$ with nodes
$\bfx_1$, $\bfx_2$, and $\bfx_3$ in the mesh
at which the angle $\angle\bfx_1\bfx_3\bfx_2$ is close to
$\pi$ (Figure~\ref{fig8} (left)).  The interpolation
$\Pi^{**} u$ for $u \in W^{2,\infty}(\tau)$ is then defined by
\begin{align*}
   (\Pi^{**} u)(\bfx_1) = u(\bfx_1), \quad
   (\Pi^{**} u)(\bfx_2) = u(\bfx_2), \quad
   (\nabla (\Pi^{**} u)(\bfx_0), \bfv) = 
   (\nabla u(\bfx_0), \bfv),
\end{align*}
where $\bfx_0$ is the foot on the edge $\overline{\bfx_1\bfx_2}$
from $\bfx_3$, and $\bfv := (\bfx_3 - \bfx_0)/|\bfx_3 - \bfx_0|$
(Figure~\ref{fig8} (right)).

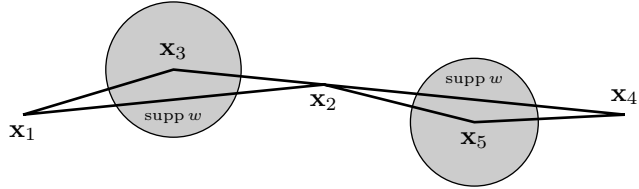
\begin{figure}[htbp]
\centering
\begin{tikzpicture}[line width = 1.0pt,scale=1.0]
   \coordinate (C) at (2.0,0.6);
   \fill [black!20] (C) circle[radius = 0.9];
   \coordinate [label=below:$\bfx_1$] (A) at (0.0,0.0);
   \coordinate [label=below:$\bfx_2$] (B) at (4.0, 0.4);
   \draw  (A) -- (B) -- (C) -- (A);
   \draw[line width = 0.5] (C) circle[radius = 0.9];
   \coordinate [label=$\bfx_3$](D) at (2.0,0.6);
   \coordinate [label={\tiny$\mathrm{supp}\,w$}](E) at (2.0,-0.25);
   \coordinate (Y) at (6.0, -0.1);
   \fill [black!20] (Y) circle[radius = 0.85];
   \draw[line width = 0.5] (Y) circle[radius = 0.85];
   \coordinate [label=above:$\bfx_4$] (X) at (8.0,0.0);
   \draw  (B) -- (X) -- (Y) -- (B);
   \coordinate [label=below:$\bfx_5$](Z) at (6.0,-0.1);
   \coordinate [label={\tiny$\mathrm{supp}\,w$}](T) at (6.0,0.25);
 \end{tikzpicture}
\caption{The supports of the correction functions should
  be disjount.}
 \label{fig9}
\end{figure}

Let $\tau_1$ be another element that shares the node $\bfx_3$.
Obviously, $(\Pi_{\tau_1}^1u)(\bfx_3) \neq (\Pi^{**}u)(\bfx_3)$ in
general.  To adjust the difference, the \textit{correction function}
$w$ is introduced so that 
\begin{gather*}
  (u + w)(\bfx_3) = (\Pi^{**}u)(\bfx_3), \qquad
  \mathrm{supp}\,w = B(\bfx_3,r_\tau), \\
  B(\bfx_3,r_\tau) := \{\bfx : |\bfx- \bfx_3| \le r_{\tau}\},
  \quad
  r_\tau := \frac{1}{2}\min\{|\bfx_3-\bfx_1|, |\bfx_3-\bfx_2|\}.
\end{gather*}
Then, $\Pi^{**}u$ is defined by $\Pi_{\tau_1}^1(u+w)$ on $\tau_1$.
To make $\Pi^{**}u$ well-defined globally, 
the supports of the correction functions should be
disjoint (see Figure~\ref{fig9}).

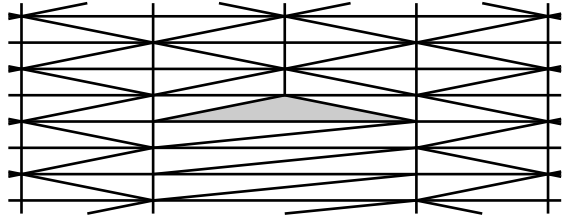
\begin{figure}[htbp]
\centering
\begin{tikzpicture}[line width = 1.0pt,scale=3.5]
   \coordinate (A1) at (0.0,0.0);
   \coordinate (A2) at (0.0,0.1);
   \coordinate (A3) at (0.0,0.2);
   \coordinate (B1) at (1.0,0.0);
   \coordinate (B2) at (1.0,0.1);
   \coordinate (B3) at (1.0,0.2);
   \coordinate (C1) at (0.5,0.3);
   \coordinate (C2) at (0.5,0.6);
   \fill [black!20] (A3) -- (C1) -- (B3) -- (A3) -- cycle;
   \draw (-0.5,-0.15) -- (-0.5,0.65);
   \draw (0.0,-0.15) -- (0.0,0.65);
   \draw (0.5,0.3) -- (0.5,0.65);
   \draw (1.0,-0.15) -- (1.0,0.65);
   \draw (1.5,-0.15) -- (1.5,0.65);
   \draw (-0.55,-0.1) -- (1.55,-0.1);
   \draw (-0.55,0.0) -- (1.55,0.0);
   \draw (-0.55,0.1) -- (1.55,0.1);
   \draw (-0.55,0.2) -- (1.55,0.2);
   \draw (-0.55,0.3) -- (1.55,0.3);
   \draw (-0.55,0.4) -- (1.55,0.4);
   \draw (-0.55,0.5) -- (1.55,0.5);
   \draw (-0.55,0.6) -- (1.55,0.6);
   \draw (-0.55,-0.01) -- (A2) -- (-0.55,0.21);
   \draw (0.0,0.3) -- (-0.55,0.41);
   \draw (-0.55,0.01) -- (0.0,-0.1) -- (-0.25,-0.15);
   \draw (-0.55,0.59) -- (-0.25,0.65);
   \draw (-0.55,0.19) -- (1.55,0.61);
   \draw (-0.55,0.39) -- (0.75,0.65);
   \draw (1.55,0.19) -- (-0.55,0.61);
   \draw (1.55,0.39) -- (0.25,0.65);
   \draw (1.55,-0.01) -- (B2) -- (1.55,0.21);
   \draw (1.0,0.3) -- (1.55,0.41);
   \draw (1.55,0.01) -- (1.0,-0.1) -- (1.25,-0.15);
   \draw (1.55,0.59) -- (1.25,0.65);
   \draw (0.0,-0.1) -- (1.0,0.0);
   \draw (1.0,-0.1) -- (0.5,-0.15);
   \draw (A1) -- (B2);
   \draw (A2) -- (B3);
   \draw (A3) -- (C1) -- (B3);
 \end{tikzpicture}
\caption{A bad element is surrounded by elements that satisfy
the maximum angle condition.}
 \label{fig10}
\end{figure}

An advantage of the approach adopted in \cite{Kucera} is that
the surrounding elements of a bad element need to satisfy only
the maximum angle condition, whereas Theorem~\ref{thm5} requires 
that such surrounding elements satisfy the minimum angle
condition (see Assumption~\ref{A2}(3)).   Hence,
\cite[Theorem~5]{Kucera} can deal with the mesh
depicted in Figure~\ref{fig10},
whereas Theorem~\ref{thm5} is not applied in this case.

In contrast, the condition that the supports of the correction
functions are disjoint is rather restrictive. For example,
\cite[Theorem~5]{Kucera} cannot deal with the bad elements
depicted in Figure~\ref{fig6} (center), whereas
Theorem~\ref{thm5} can be applied.
In \cite[Section~3.4]{Kucera}, Ku\v{c}era presented a way to
deal with clusters of bad elements.
However, each such cluster must be sufficiently far
from other clusters and the boundary $\partial\Omega$.
Note that Theorem~\ref{thm5} can be applied in the cases
of a clusters of bad elements that are close or connected to
other clusters or the boundary $\partial\Omega$, if
Assumption~\ref{A2}(3) is satisfied.

Moreover, the approach in \cite{Kucera} deals with only 
triangular elements ($n=2$) in the $W^{2,\infty}$
setting, whereas Theorem~\ref{thm5} is applicable to the more general 
case of $n$-simplices ($n=2,3$) in the $H^2$ setting. 

Second, we explain the approach of Duprez, Lleras, and Lozinski
\cite{DLL}, which is slightly
different from that of \cite{Kucera}.
Their assumption is that each cluster of simplices that contains
bad elements should be star-shaped and completely surrounded by good simplices.

\begin{figure}[htbp]
\centering
\begin{tikzpicture}[line width = 1.0pt,scale=0.7]
 %  \draw [help lines] (-2,-2) grid (3,3);
   \coordinate (A) at (-1.7,-1.3);
   \coordinate (B) at (0.1,-1.8);
   \coordinate (C) at (-0.2,-0.1); 
   \coordinate (D) at (1.3,0.0);
   \coordinate (E) at (1.9,-1.9);
   \coordinate (F) at (2.1,0.7);
   \coordinate (G) at (3.2,-0.9);
   \coordinate (H) at (3.3,0.8);
   \coordinate (I) at (2.2,2.3);
   \coordinate (J) at (0.1,1.6);
   \coordinate (K) at (0.2,3.0);
   \coordinate (L) at (-1.8,2.0);
   \coordinate (M) at (-2.0,0.5);
   \fill [black!20] (C) -- (D) -- (F) -- (C) -- cycle;
   \fill [black!8] (C) -- (F) -- (J) -- (C) -- cycle;
   \draw (A) -- (B) -- (C) -- (A);
   \draw (C) -- (D) -- (B);
   \draw (B) -- (E) -- (D);
   \draw (C) -- (F) -- (D);
   \draw (D) -- (G) -- (E);
   \draw (F) -- (G) -- (H) -- (F);
   \draw (H) -- (I) -- (F);
   \draw (C) -- (J) -- (F);
   \draw (I) -- (J) -- (K) -- (I);
   \draw (K) -- (L) -- (J) -- (M) -- (C);
   \draw (A) -- (M) -- (L);
\end{tikzpicture}
\caption{A bed element is surrounded by good elements. }
\end{figure}
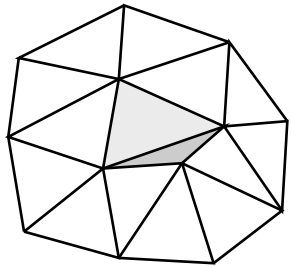

The interpolation $\Pi^{***}u$ is then defined by Taylor polynomials on
the cluster.  On the surrounding good elements, the Lagrange
interpolation $\Pi_{\tau}^1u$ is modified to ensure the
continuity of $\Pi^{***}u$ between elements. To evaluate the influence
of the procedure, Duprez, Lleras, and Lozinski used $L^\infty$ error
estimates of the Sobolev embedding, which is quite similar to our
approach.  Their main theorem \cite[Proposition~2.7]{DLL} can deal with
both two- and three-dimensional cases in the $H^2$ setting.

Although \cite[Proposition~2.7]{DLL} is applicable to a wide range
of configurations of bad elements, it seems that Theorem~\ref{thm5}
can also deal with most of such cases.  Furthermore, as
Assumption~\ref{A2} does not require that each cluster of bad elements
is surrounded by good elements, Theorem~\ref{thm5} is applicable to
a wider range of configurations than \cite[Proposition~2.7]{DLL}.
For example, the mesh given in the numerical experiment 
(Figure~\ref{fig30}) cannot be dealt with using the cited approach.

\section{Numerical experiment}
We present a numerical experiment to confirm the theoretical results
obtained in the previous sections. Let $\Omega := (0,1) \times (0,1)$.
We consider the Poisson equation
\begin{align}
   - \Delta u = 2(x(1-x) + y(1-y)) \;\text{ in } \; \Omega, \qquad
    u = 0 \;\text{ on }\; \partial\Omega,
 \label{num-exam}
\end{align}
whose exact solution is $u = x(1-x)y(1-y)$.

\begin{figure}[htbp]
\centering
\begin{tikzpicture}[line width = 0.8pt,scale=0.4]
\coordinate (X) at (2.0,0.0);
\coordinate (Y) at (2.0,2.0);
\coordinate (Z) at (0.0,2.0);
\coordinate (S) at (1.0,0.2);
\coordinate (T) at (1.0,1.8);
\foreach \y in {0, 2, 4, 6, 8}
{
\foreach \x in {0, 2, 4, 6, 8}
{\coordinate (A) at (\x,\y);
  \coordinate (B) at ($(A)+(X)$);
  \coordinate (C) at ($(A)+(Y)$);
  \coordinate (D) at ($(A)+(Z)$);
  \coordinate (E) at ($(A)+(S)$);
  \coordinate (F) at ($(A)+(T)$);
  \draw (B) -- (C) -- (E) -- (B) -- (A) -- (E) -- (F) -- (C) -- (D) --
 (F) -- (A) -- (D);
  };};
\end{tikzpicture} 
\quad
\begin{tikzpicture}[line width = 1.0pt,scale=0.4]
   \coordinate [label=below:{\small$(x,y)$}](A) at (0.0,0.0);
   \coordinate [label=below:{\small$(x+k,y)$}](B) at (7.0,0.0);
   \coordinate [label=above:{\small$(x+k,y+k)$}](C) at (7.0,7.0);
   \coordinate [label=above:{\small$(x,y+k)$}](D) at (0.0,7.0);
   \coordinate (E) at (3.5,0.8);
   \coordinate (F) at (3.5,6.2);
  \draw (B) -- (C) -- (E) -- (B) -- (A) -- (E) -- (F) -- (C) -- (D) --
 (F) -- (A) -- (D);
   \coordinate [label=right:{\small$(x+\frac{k}{2},y+\alpha k)$}](G) at (7.3,0.9);
   \draw [dotted, line width = 0.7pt][->] (G) -- (3.8,0.9);
   \coordinate [label=right:{\small$(x+\frac{k}{2},y+(1-\alpha) k)$}](H) at (7.3,6.2);
   \draw [dotted, line width = 0.7pt][->] (H) -- (3.8,6.1);
\end{tikzpicture} 
\caption{Triangulation $\T_h$ of $\Omega$ (left), and
the subdivision of each small square (right). }
\label{fig30}
\end{figure}
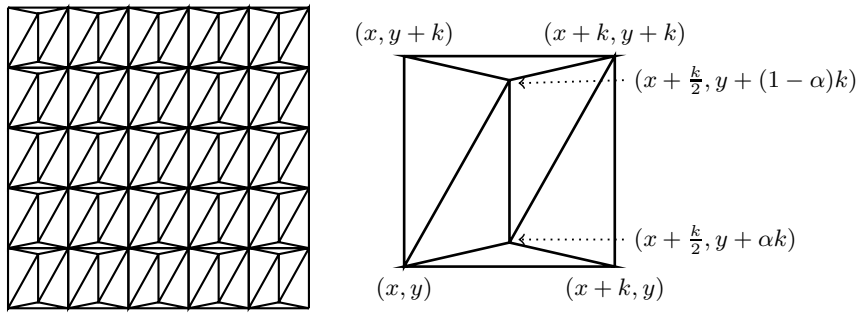

Let $K$ be a positive integer and $k = 1/K$.  To define
the triangulation $\T_h$, we {first} divide $\Omega$ into $K^2$ small squares
with side lengths $k$.  Each small square is further {sub}divided into {six}
triangles as depicted in Figure~\ref{fig30}.   Here, the coordinates of the
points in each small square are $(x,y)$, $(x+k,y)$, $(x+k,y+k)$,
$(x,y+k)$, {$(x+k/2,y + \alpha k)$, and $(x+k/2,y + (1 - \alpha) k)$}, where
$0 < \alpha \ll 1$.  Note that the smallest triangle becomes very bad when
$\alpha$ approaches to zero.  It is clear that this triangulation
satisfies Assumption~\ref{A5}, and Corollary~\ref{cor7} is applicable to
this triangulation $\T_h$.

We compute the piecewise linear finite element solutions
$u_h \in S_{h0}$ of \eqref{num-exam} on $\T_h$, and the error
$|u - u_h|_{H^1(\Omega)}$ with varying $K$ and $\alpha$.
Note that $h = (0.25 + (1-\alpha)^2)^{1/2}k$ if
$\alpha < (2-\sqrt{3})/2 = 0.1339\dots$.

Table~1 presents the results of the numerical experiment.
Here, $e_h$ stands for $|u - u_h|_{H^1(\Omega)}$.
Although the meshes $\T_h$ contain many very bad elements,
we observe that the error $|u - u_h|_{H^1(\Omega)}$ is of
$\OO(h)$ and consistent with Theorem~\ref{thm5} and Corollary~\ref{cor7}.

\vspace{3mm}
\begin{center}
\begin{table}[htpb]
\caption{Results of the numerical experiment.  In the table,
$e_h$ stands for $|u - u_h|_{H^1(\Omega)}$.}
\begin{tabular}{|c|c|c|c|c|c|c|} \hline
  & \multicolumn{2}{|c|}{$\alpha=0.1$} 
  & \multicolumn{2}{|c|}{$\alpha=0.01$}
  & \multicolumn{2}{|c|}{$\alpha=0.0001$} \\ \hline
 $K$ & $e_h$ & $e_h/h$ &  $e_h$ & $e_h/h$ & $e_h$ & $e_h/h$ \\ \hline
 10 & 1.8002 e-2 & 0.17485 & 2.0839 e-2 & 0.18789 & 2.1237 e-2 & 0.18997\\ \hline
 20 & 9.0151 e-3 & 0.17512 & 1.0440 e-2 & 0.18827 & 1.0641 e-2 & 0.19036\\ \hline
 40 & 4.5093 e-3 & 0.17519 & 5.2229 e-3 & 0.18836 & 5.3231 e-3 & 0.19046\\ \hline
 80 & 2.2548 e-3 & 0.17521 & 2.6118 e-3 & 0.18839 & 2.6619 e-3 & 0.19049\\ \hline
160 & 1.1274 e-3 & 0.17521 & 1.3059 e-3 & 0.18866 & 1.3310 e-3 & 0.19049\\ \hline
\end{tabular}
\end{table}
\end{center}

\section{Concluding remarks and future works}
In this paper, we examined the error estimation
of finite element solutions on meshes that contain many
bad elements.  Our research was motivated by the works of
{Ku\v{c}era \cite{Kucera} and}
Duprez, Lleras, and Lozinski \cite{DLL}.
We showed that even if $\T_h$ contains many
bad elements, finite element solutions $\{u_h\}$ converge
to the exact solution $u \in H^2(\Omega)$ with order
$\OO(h)$, if each bad element is covered virtually by a
simplex {that satisfies the minimum angle condition}.
  To be precise, we presented Assumption~\ref{A2}
that is very general and covers a wide range of triangulation patterns.

For the case $n = 2$, we gave Assumption~\ref{A5}, which
is simpler than Assumption~\ref{A2} but is still practical
and applicable to many meshes that cannot be dealt with
using the methods of {\cite{Kucera, DLL}}. 

In our numerical experiment, the finite element solutions
behaved exactly as predicted from Theorem~\ref{thm5} and
Corollary~\ref{cor7}.

The results of this paper can be extended in many ways.  An immediate
extension may be the consideration of finite element solutions using
higher order elements, and the authors are now working on this
extension.

%In the numerical example given
%in Section~5, we see that the ratio of the areas
%\begin{align*}
%   \frac{|\Omega^{bad}(\T_h)|}{|\Omega|}, \qquad
%    \Omega^{bad}(\T_h) := \bigcup_{\tau \in \T_h
%    \text{ is a bad element}} \tau
%\end{align*}
%decreases as $\alpha \to 0$. 
%This observation leads the following conjecture. 
%Suppose that we have a sequence $\{\T^n\}_{n=1}^\infty$ of meshes that
%satisfies
%\begin{align*}
%  h_n := \max_{\tau \in \T^n} h_\tau, \qquad  
%  \lim_{n \to \infty} \frac{|\Omega^{bad}(\T^n)|}{|\Omega|} = 0.
%\end{align*}
%We might then be able to show that finite element solutions
%$u_n$ defined on each $\T^n$ still enjoy the standard error estimation
%\begin{align*}
%  |u - u_n|_{H^1(\Omega)} \le C h_n |u|_{H^2(\Omega)},
%\end{align*}
%where $C$ is a constant independent of $n$.

{\small
{\em Authors' addresses}: \\
{\em Kenta Kobayashi}, Hitotsubashi University, Kunitachi, Japan \\
 e-mail: \texttt{kenta.k@\allowbreak r.hit-u.ac.jp}. \\
{\em Takuya Tsuchiya}, Osaka University, Toyonaka, Japan \\
 e-mail: \texttt{tsuchiya.takuya.plateau@kyudai.jp}.
}

\end{document}